\documentclass[a4paper]{amsart}
\usepackage{t-angles}
\usepackage{xypic}
\newcommand{\catC}{\mathcal{C}}
\newcommand{\catB}{\mathcal{B}}
\newcommand{\Z}{\mathcal{Z}}
\newcommand{\ot}{\otimes}

\newtheorem{propo}{Proposition}[section]

\theoremstyle{definition}

\theoremstyle{remark}
\newtheorem{remar}[propo]{Remark}

\DeclareMathOperator{\length}{\mathsf{length}}
\DeclareMathOperator{\End}{\mathsf{End}}
\DeclareMathOperator{\Hom}{\mathsf{Hom}}
\DeclareMathOperator{\id}{\mathsf{id}}
\newcommand{\De}{\Delta}
\newcommand{\fect}[1]{{\fontencoding{OT1}\fontfamily{pcr}\fontshape{it}\fontseries{b}\fontsize{12}{0}\selectfont%
	\text{#1\/}}}

\def\xcirc{\objectmargin{0.1pc}%
	\objectstyle{\sssize}\diagram\squarify<1pt>{}\circled\enddiagram}

\begin{document}
\title{How to disentangle two braided Hopf algebras}
\author{M. Gra\~na}
\author{J.A. Guccione}
\author{J.J. Guccione}
\address{%
M.G., J.A.G., J.J.G.:\newline\indent
    Depto de Matem\'atica - FCEyN\newline\indent
    Universidad de Buenos Aires\newline\indent
    Pab. I - Ciudad Universitaria\newline\indent
    1428 - Buenos Aires - Argentina\newline
	\email{M.G.}{matiasg@dm.uba.ar}
	\email{J.A.G.}{vander@dm.uba.ar}
	\email{J.J.G.}{jjgucci@dm.uba.ar}}
\thanks{This work was partially supported by CONICET, ANPCyT (PICT-02 12330 and PICT-04 20054), UBA X294}
\begin{abstract}
	We show how to define the tensor product of two braided Hopf algebras.
\end{abstract}
\maketitle

\section*{Introduction}
In this short note we show how to endow with a canonical Hopf algebra structure the tensor product
$H\otimes L$ of two braided Hopf algebras living in a monoidal category, provided that there is an
isomorphism $c:L\otimes H\to H\otimes L$ satisfying suitable conditions (see
Proposition~\ref{prop1.1}). In particular, the square tensor $H\otimes H$ is always a braided Hopf
algebra with this structure. Also, the tensor product of two Hopf algebras living in a braided
category is again a braided Hopf algebra (although in a new braided category). When seen in this
context, the key point is to replace the braid $c_{HL}^{\phantom{-1}}$ with $c_{LH}^{-1}$.

\section{Preliminaries}
We work in a monoidal category $\catC$, for instance, the category of vector spaces over a field
$k$. We write $\otimes$ and  $I$ for the tensor product and the unit of $\catC$, respectively. The
associativity and unit constraints are assumed without referring to them. We assume the reader is
familiar with the notions of algebras and coalgebras in monoidal categories.  All the algebras are
associative unitary and the coalgebras are coassociative counitary. Given an algebra $A$ and a
coalgebra $C$, we let $\mu\colon A\ot A \to A$, $\eta\colon I \to A$, $\De\colon C\to C\ot C$ and
$\varepsilon\colon C\to I$ denote the multiplication, the unit, the comultiplication and the counit,
respectively, specified with a subscript if necessary. We are going to use the nowadays well known
graphic calculus for monoidal and braided categories. As usual, morphisms will be composed downwards
and tensor products will be represented by horizontal concatenation in the corresponding order. The
identity map of an object will be represented by a vertical line.  Given an algebra $A$, the
diagrams
$\begin{tangle}\HH\hcu\end{tangle}$
	and
\raisebox{-2pt}{$\begin{tangle}\HH\unit\end{tangle}$}
stand for the multiplication map and the unit of $A$ respectively. Given a
coalgebra $C$, the comultiplication and the counit of $C$ will be represented
by the diagrams
$\begin{tangle}\HH\hcd\end{tangle}$
	and
\raisebox{-8pt}{$\begin{tangle}\HH\counit\end{tangle}$}
respectively.

Let $V$, $W$ be objects in $\mathcal{C}$ and let $c\colon V\ot W \to W\ot V$ be an arrow.
\begin{itemize}
\item If $V$ is an algebra, then we say that $c$ is {\em compatible} with the
	algebra structure of $V$ if $c \xcirc (\eta\ot W)= W\ot \eta$ and
	$c\xcirc (\mu\ot W)=  (W\ot \mu)\xcirc(c\ot V)\xcirc(V\ot c)$.
\item If $V$ is a coalgebra, then we say that $c$ is compatible with the
	coalgebra structure of $V$ if $(W\ot\varepsilon)\xcirc c=\varepsilon\ot W$
	and $(W\ot \De) \xcirc c = (c\ot V)\xcirc (V\ot c)\xcirc (\De \ot W)$.
\end{itemize}
\noindent Of course, there are similar compatibilities when $W$ is an algebra or a coalgebra.

\section{Tensor products of braided Hopf algebras}
Recall that a braided bialgebra in $\catC$ is an object $H$ of $\catC$ endowed with an algebra
structure, a coalgebra structure and an isomorphism $c_H\in \End_{\catC} (H^2)$ (called the braid of $H$)
satisfying the Braid Equation
\[
(c_H\otimes H)\xcirc (H\otimes c_H)\xcirc (c_H\otimes H)
	=(H\otimes c_H)\xcirc (c_H\otimes H)\xcirc (H\otimes c_H),
\]
such that: $c_H$ is compatible with the algebra and coalgebra structures of $H$, $\eta$ is a
coalgebra morphism, $\varepsilon$ is an algebra morphism and
$\Delta\xcirc\mu = (\mu\ot \mu)\xcirc(H\ot c_H \ot H)\xcirc(\Delta \ot \Delta)$.
If moreover there exists a map $S\colon H\to H$ which is the inverse of the identity in the monoid
$\Hom_{\catC}(H,H)$ with the convolution product, then we say that $H$ is a braided Hopf algebra and
we call $S$ the antipode of $H$.

\smallskip

Let $H$ be a braided bialgebra in $\mathcal{C}$. It is well known that if the braid of $H$ is
involutive (i.e., $c_H^2=\id$), then $H\ot H$ is a braided bialgebra in a natural way. The following
Proposition is the main result in this note. It shows in particular that the involutivity hypothesis
can be removed.

\begin{propo}\label{prop1.1}
	Let $H$ and $L$ braided bialgebras in $\mathcal{C}$ and let
	$c_{LH}\colon L\ot H\to H\ot L$ be a invertible arrow. If
	\begin{align*}
		& (c_{LH}\ot L)\xcirc (L\ot c_{LH}) \xcirc (c_L\ot H)
			=(H\ot c_L)\xcirc (c_{LH}\ot L)\xcirc (L\ot c_{LH}),\\
		& (c_H\ot L)\xcirc (H\ot c_{LH})\xcirc (c_{LH}\ot H)
			= (H\ot c_{LH})\xcirc(c_{LH}\ot H)\xcirc (L\ot c_H),
	\end{align*}
	and $c_{LH}$ is compatible with the bialgebras structures of $H$ and $L$,
	then $H\ot L$ is a braided bialgebra, via
	\begin{itemize}
		\item $\mu_{H\ot L}= (\mu_H \ot \mu_L)\xcirc (H\ot c_{LH}\ot L),$
		\item $\eta_{H\ot L}=\eta_H\ot \eta_L$,
		\item $\Delta_{H\ot L}= (H\ot c_{LH}^{-1}\ot H)\xcirc (\Delta_H\ot \Delta_L)$,
		\item $\varepsilon_{H\ot L}= \varepsilon_H\ot\varepsilon_L$,
		\item $c_{H\ot L}= (H\ot c_{LH}^{-1}\ot L)\xcirc (c_H\ot c_L)\xcirc (H\ot c_{LH}\ot L)$.
	\end{itemize}
	Moreover, if $H$ and $L$ are braided Hopf algebras, then so is $H\ot L$,
	with antipode $S_{H\ot L}=S_H\ot S_L$.
\end{propo}
\begin{proof}
	The fact that $H\otimes L$ is an algebra and a coalgebra is well-known and
	standard.
	We check now the compatibility between multiplication and
	comultiplication:
	
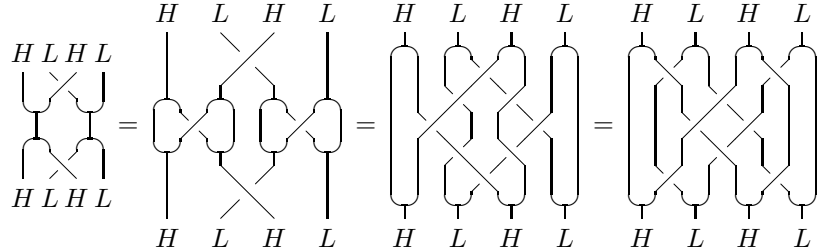
\begin{align*}
		\begin{tangles}{lcr}
			\HH\object{H} & \object{L} \step[1]\object{H} & \step[1]\object{L} \\
			\HH\id & \xx & \step[1]\id \\
			\HH \cu && \cu \\
			\HH \cd && \cd \\
			\HH\id & \x & \step[1]\id \\
			\HH\object{H} & \object{L} \step[1]\object{H} & \step[1]\object{L} \\
		\end{tangles}
		\;=\;
		\begin{tangles}{lccr}
			\HH\step[.5]\object{H} \step[2] \object{L} \step[2]\object{H} \step[2]\object{L} \\
			\step[.5]\id\step[2]\xx\step[2]\id \\
			\HH\cd \step[1]\cd \step[1] \cd \step[1] \cd \\
			\HH \id\step[1] \xx \step[1] \id\step[1]  \id\step[1]  \xx \step[1] \id \\
			\HH\cu \step[1]\cu \step[1] \cu \step[1] \cu \\
			\step[.5]\id\step[2]\x\step[2]\id \\
			\HH\step[.5]\object{H} \step[2] \object{L} \step[2]\object{H} \step[2]\object{L} \\
		\end{tangles}
		\;=\;
		\begin{tangles}{lccr}
			\HH\step[.5]\object{H} \step[2] \object{L} \step[2]\object{H} \step[2]\object{L} \\
			\HH\cd \step[1]\cd \step[1] \cd \step[1] \cd \\
			\HH\id\step[1]\id\step[1]\id\step[1]\xx\step[1]\id\step[1]\id\step[1]\id \\
			\HH \id\step[1]  \id\step[1]  \xx \step[1] \xx \step[1] \id\step[1]  \id \\
			\HH \id\step[1] \xx \step[1] \id\step[1]  \id\step[1]  \xx \step[1] \id \\
			\HH \id\step[1]\id\step[1]\x\step[1]\x\step[1]\id\step[1] \id \\
			\HH\id\step[1]\id\step[1]\id\step[1]\x\step[1]\id\step[1]\id\step[1]\id \\
			\HH\cu \step[1]\cu \step[1] \cu \step[1] \cu \\
			\HH\step[.5]\object{H} \step[2] \object{L} \step[2]\object{H} \step[2]\object{L} \\
		\end{tangles}
		\;=\;
		\begin{tangles}{lccr}
			\HH\step[.5]\object{H} \step[2] \object{L} \step[2]\object{H} \step[2]\object{L} \\
			\HH\cd \step[1]\cd \step[1] \cd \step[1] \cd \\
			\HH\id \step[1]\x\step[1]\id\step[1]\id\step[1]\x\step[1]\id  \\
			\HH\id\step[1]\id\step[1]\id\step[1]\xx\step[1]\id\step[1]\id\step[1]\id \\
			\HH\id\step[1]\id\step[1]\xx\step[1]\xx\step[1]\id\step[1]\id \\
			\HH\id\step[1]\id\step[1]\id\step[1]\x\step[1]\id\step[1]\id\step[1]\id \\
			\HH\id \step[1]\xx \step[1]\id \step[1] \id \step[1]\xx \step[1]\id  \\
			\HH\cu \step[1]\cu \step[1] \cu \step[1] \cu \\
			\HH\step[.5]\object{H} \step[2] \object{L} \step[2]\object{H} \step[2]\object{L} \\
		\end{tangles}
	\end{align*}
	We leave to the reader to prove that $c_{H\ot L}$ satisfies the Braid
	Equation, and that it is compatible with $\mu_{H\ot L}$, $\Delta_{H\ot L}$,
	$\eta_{H\ot L}$ and $\varepsilon_{H\ot L}$.
	The proof of the last assertion in the statement is also straightforward.
\end{proof}

\section{Braided families and compatible maps}
The aim of this section is to give a more categorical proof of
Proposition~\ref{prop1.1}.  The methods presented here could be useful in
generalizing this result to braided versions of bicrossproducts, matched pairs,
etc. (see \cite{MR1706813}).  Let $\fect V=(V_i)_{i\in\Im}$ be a family of objects in
$\catC$. A family  $\fect C$ of isomorphisms
$$
c_{ij}:V_i\otimes V_j\to V_j\otimes V_i\quad \text{($i,j\in\Im$)}
$$
is said to be {\em braided} if $\forall i,j,k\in\Im$ the Braid Equations
\begin{equation*}
    (V_k\otimes c_{ij})\xcirc (c_{ik}\otimes V_j)\xcirc (V_i\otimes c_{jk})
        =(c_{jk}\otimes V_i)\xcirc (V_j\otimes c_{ik})\xcirc (c_{ij}\otimes V_k).
\end{equation*}
are satisfied. Given such a family and a string $\mathbf{i}=(i_1,\dots, i_n)$ of elements in
$\Im$, we call $n$ the {\em length} of $\mathbf{i}$, we let $\mathbf{i}_{<n}$ denote the string
$(i_1,\dots, i_{n-1})$ and we put $V_{\mathbf{i}}=V_{i_1}\otimes \cdots \otimes V_{i_n}$. For each
pair $\mathbf{i},\mathbf{j}$ of such strings we define the map $c_{\mathbf{i}\mathbf{j}}\colon
V_{\mathbf{i}}\otimes V_{\mathbf{j}}\rightarrow V_{\mathbf{j}}\otimes V_{\mathbf{i}}$, recursively
by

\begin{itemize}
\setlength{\itemsep}{8pt}
\item If $\mathbf{i}=(i)$ and $\mathbf{j}=(j)$, then $c_{\mathbf{ij}}=c_{ij}$.
\item If $\mathbf{i}=(i)$ and $\length(\mathbf{j})=m>1$, then
	$c_{\mathbf{ij}}=(V_{\mathbf{j}_{<m}} \otimes c_{ij_m})
		\xcirc(c_{\mathbf{ij}_{<m}} \otimes V_{j_m})$.
\item If $\length(\mathbf{i})=n>1$, then
	$c_{\mathbf{ij}}=(c_{\mathbf{i}_{<n},\mathbf{j}}\otimes V_{i_n})
		\xcirc (V_{\mathbf{i}_{<n}}\otimes c_{(i_n)\mathbf{j}})$.
\end{itemize}
We write $c_{i\mathbf{j}}=c_{(i)\mathbf{j}}$ and $c_{\mathbf{i}j}=c_{\mathbf{i}(j)}$.

We say that a map $f\colon V_{\mathbf{i}}\rightarrow V_{\mathbf{j}}$ is
\emph{compatible} with $\fect C$ if
\[
    (f\otimes V_l)\xcirc c_{\mathbf{i}l} = c_{\mathbf{j}l}\xcirc (V_l\otimes f)
		\quad\text{and}\quad
    (V_l \otimes f)\xcirc c_{l\mathbf{i}} = c_{l\mathbf{j}}\xcirc ( f \otimes V_l)
\]
for each $l\in \Im$.

Let $\mathfrak{D}=( \fect V, \fect C,\fect M)$, where $\fect M$ is a family of
maps compatible with $\fect C$. We want to embed this datum into a braided category in a
natural way. To this end, let $\catB_{\mathfrak{D}}$ be the category whose objects are pairs
$\mathbf{W}=(W,\lambda^{W}_-)$, where $W$ is an object in $\catC$ and
$\lambda^{W}_-=(\lambda^{W}_i)_{i\in\Im}$ is a family  of isomorphisms
$\lambda^{W}_i:V_i\otimes W\to W\otimes V_i$, subject to
\begin{itemize}
	\item $(W\otimes c_{ij})\xcirc \lambda^{W}_{(ij)}
		= \lambda^{W}_{(ji)}\xcirc (c_{ij}\otimes W)$,
	\item $(W\otimes f)\xcirc \lambda_{\mathbf{i}}^W=\lambda_{\mathbf{j}}^W\xcirc (f\otimes W)$
		for each map $f\colon V_{\mathbf{i}}\to V_{\mathbf{j}}$ in $\fect M$,
\end{itemize}
where $\lambda^W_{\mathbf{i}}\colon V_{\mathbf{i}}\otimes W\to W\otimes V_{\mathbf{i}}$
is recursively defined as follows:
\begin{itemize}
	\item If $\mathbf{i}=(i)$, then $\lambda^W_{\mathbf{i}}=\lambda^W_i$.
	\item If $n=\length(\mathbf{i})>1$, then
		$\lambda^W_{\mathbf{i}}=(\lambda^W_{\mathbf{i}_{<n}} \otimes V_{i_n})
			\xcirc (V_{\mathbf{i}_{<n}}\otimes \lambda^W_{i_n})$
\end{itemize}

A morphism $\phi:\mathbf{W}\to\mathbf{Z}$ in $\catB_{\mathfrak{D}}$ is a map $\phi:W\to Z$ in
$\catC$ such that
\[
(\phi\otimes V_i)\xcirc \lambda^{W}_i
    =\lambda^{Z}_i\xcirc (V_i\otimes\phi)\quad\text{for all $i\in\Im$.}
\]

The category $\catB_{\mathfrak{D}}$ is monoidal, with
\begin{itemize}
    \item tensor product given by
        $\mathbf{W}\otimes\mathbf{Z} =(W\otimes Z,(W\otimes\lambda^{Z}_-)\xcirc
            (\lambda^{W}_-\otimes Z))$,
    \item unit $\mathbf{I}=(I,\lambda^{I}_-)$, where
        $\lambda^{I}_i:V_i\otimes I\to I\otimes V_i$
        is the canonical map in $\catC$.
    \item associativity and unit constraints induced from those of $\catC$.
\end{itemize}

\begin{remar}
    The family $\fect V$ can be embedded in $\catB_{\mathfrak{D}}$ by means of
	$\mathbf{V}_j=(V_j,\lambda^{{V}_j}_-)$, where ${\lambda^{{V}_j}_i=c_{ij}}$.
\end{remar}

We consider $\Z(\catB_{\mathfrak{D}})$, the center of $\catB_{\mathfrak{D}}$, which, we recall
from \cite[Cor. 9.1.6]{MR1381692}, has as objects the pairs $(\mathbf{W},\gamma_{\mathbf{W},-})$,
where $\mathbf{W}$ is an object in $\catB_{\mathfrak{D}}$ and $\gamma_{\mathbf{W},-}$ is a
natural isomorphism $(\mathbf{W}\otimes -)\longrightarrow (-\otimes \mathbf{W})$ such that
\begin{align*}
    &\gamma_{\mathbf{W},\mathbf{I}} \quad\text{is the canonical isomorphism in
    $\catB_{\mathfrak{D}}$,} \\
    &\gamma_{\mathbf{W},\mathbf{X}\otimes \mathbf{Y}}
        =(\mathbf{X}\otimes\gamma_{\mathbf{W},\mathbf{Y}})\xcirc
            (\gamma_{\mathbf{W},\mathbf{X}}\otimes \mathbf{Y})
				\quad\forall\mathbf{X},\mathbf{Y}\in\catB_{\mathfrak{D}}.
\end{align*}
An arrow $\boldsymbol{\theta}\colon (\mathbf{W},\gamma_{\mathbf{W},-})\to (\mathbf{Z},\gamma_{\mathbf{Z},-})$
is an arrow $\theta\colon \mathbf{W}\to\mathbf{Z}$ in $\catB_{\mathfrak{D}}$ such that
$$
	(\mathbf{X}\otimes\theta)\xcirc\gamma_{\mathbf{W},\mathbf{X}}
		=\gamma_{\mathbf{Z},\mathbf{X}}\xcirc(\theta\otimes\mathbf{X})
			\quad\forall\mathbf{X}\in\catB_{\mathfrak{D}}.
$$

\begin{propo}
    \label{pr:32}
    The initial datum $\mathfrak{D}$ is included in the center via the identification
	$\iota\colon\mathfrak{D}\to\Z(\catB_{\mathfrak{D}})$ given by
    \begin{itemize}
    \setlength{\itemsep}{4pt}
        \item $\iota V_j=(\mathbf{V}_j,\gamma_{\mathbf{V}_j,-})$, where
            $\gamma_{\mathbf{V}_j,\mathbf{W}}=\lambda^{{W}}_j$ for all $V_j\in\fect V$.
		\item $\iota c_{ij}=\gamma_{\mathbf{V}_i,\mathbf{V}_j}=\lambda^{{V}_j}_i=c_{ij}$ for all
			$c_{ij}\in\fect C$.
		\item $\iota f=f$ for all $f\in\fect M$.
    \end{itemize}
\end{propo}
\begin{proof}
    Straightforward.
\end{proof}

\subsection{Alternative proof of Proposition~\ref{prop1.1}}
\begin{proof}
	We take $\fect V=(V_1=H,\,V_2=L,\,V_0=I)$. $\fect C$ consists of the arrows ${c_{11}=c_{H}}$,
	$c_{12}=c_{HL}$, $c_{21}=(c_{HL})^{-1}$, $c_{22}=c_{L}$ and the canonical maps
	$V_0\otimes V_i\to V_i\otimes V_0$, $V_i\otimes V_0\to V_0\otimes V_i$.
	We take also $\fect M=\{\varepsilon_H,\varepsilon_L,\eta_H,\eta_L,\Delta_H,\Delta_L,\mu_H,\mu_L\}$.
	By Proposition~\ref{pr:32}, these data are included, via a map $\iota$, in a braided category.
	Moreover, the braiding between $\iota V_1$ and $\iota V_2$ is involutive, i.e.,
	$c_{\iota V_1,\iota V_2}\xcirc c_{\iota V_2,\iota V_1}=\id$. Hence, the first item of
	\cite[Corollary~2.17]{MR1706813} applies to give the desired result.
\end{proof}

\end{document}